\documentclass[a4paper]{article}
\usepackage{latexsym}
\usepackage{amssymb}
\usepackage{amsmath}
\newtheorem{conj}{Conjecture}[section]
\newtheorem{lemma}[conj]{Lemma}
\newtheorem{coro}[conj]{Corollary} 
\newtheorem{defi}[conj]{Definition}

\newtheorem{thm}{Theorem}
\newtheorem{rem}[conj]{Remark}

\newcommand{\qed}{\raisebox{-.8ex}{$\Box$}}
\newenvironment{bew}
{\noindent{\bf Proof.}}
{\hfill \qed\\}

\newcommand{\la}{\langle}
\newcommand{\ra}{\rangle}

\newcommand{\Aut}{{\rm Aut}}

\newcommand{\Syl}{{\rm Syl}}
\newcommand{\PGL}{{\rm PGL}}
\newcommand{\PSL}{{\rm PSL}}

\newcommand{\Stab} {{\rm Stab}}

\newcommand{\Sym}{~{\rm Sym}}

\newcommand{\Alt}{~{\rm Alt}}
\newcommand{\RMult}{~{\rm RMult}}

\newcommand{\NN}{\mathbb{N}}
\newcommand{\QR}{{\rm I}\kern-5.0pt {\rm Q} \kern2pt}

\title{The finite Bruck Loops\thanks{This research is part of the project ``Transversals in Groups with an application to loops'' GZ: BA 2200/2-2 funded by
the DFG}}
\author{B. Baumeister, A. Stein}
\begin{document}

\maketitle

\begin{abstract}
We continue the work by Aschbacher, Kinyon and Phillips [AKP] as well as of Glauberman [Glaub1,2]
by describing the structure of the finite Bruck loops. We show that a finite Bruck loop $X$ 
is the direct product of a Bruck loop of odd order with either a soluble Bruck loop of $2$-power order
or a product of loops related to the groups $PSL_2(q)$, $q= 9$ or $q \geq 5$ a Fermat prime. The latter possibility
does occur as is shown in [Nag1, BS]. As corollaries we obtain versions of Sylow's, Lagrange's and Hall's
Theorems for loops.
\end{abstract}

\section{Introduction}
Let $(X,\circ)$ be a finite loop; that is a finite set together with a binary operation $\circ$ on $X$, such that
there exists an element $1 \in X$ with $1 \circ x = x \circ 1 =x$ for all $x \in X$ and such that
the left and right translations \[ \lambda_x : X \to X ,~ y \mapsto x \circ y,\quad \rho_y : X \to X,~ x \mapsto x \circ y \]
are bijections. Loops can be thought of as groups without associativity law.

Given a loop $X$, let $G:= \langle \rho_x: x \in X \rangle \le \Sym(X)$, the so called {\em  enveloping group of $X$}.
 The set $K:=\{ \rho_x: x \in X \}$ is a transversal to $H:=\Stab_G(1)$
and $(G,H,K)$ is called the {\em  Baer envelope of $X$}.
This connection between loops and  transversals in groups goes back to Baer \cite{B}, see Section 2.
In [Asch] Aschbacher started the study of loops using group theory [Asch], which turned out to be a 
very powerful tool [Asch, AKP, Nag1, BS].

Though loops are a generalisation of groups, general loops can be very wild due to the missing associativity: 
Left-and right inverses may not be identical, powers of elements may not be definable in the usual way and 
many loops without proper subloops besides the cyclic groups of prime order exist. 

In search for natural restrictions on loops, Bol discovered in \cite{Bol} the following identity, which today is known as the {\em (right) Bol identity}:
\[\quad ((x \circ y) \circ z) \circ y = x \circ ((y \circ z) \circ y) \mbox{ for all }~x,y,z \in X. \]

A loop is called a (right) {\em  Bol loop}, if it satisfies the above identity.
Bol himself showed, that this generalisation of associativity is quite natural, but that groups
are not the only examples of Bol loops. 

One consequence of the identity is, 
that the subloop generated by one element is a (cyclic) group. Therefore powers and inverses of elements are well defined. 

Examples of Nagy [Nag2] however imply, that general Bol loops may still be quite wild:
While groups of odd order are soluble due to work of Feit and Thompson, simple Bol loops of odd non-prime order exist. 
Furthermore there are noncyclic simple Bol loops, which occur from transversals in soluble groups.

A natural further restriction is the following identity known as the {\em  automorphic inverse property} AIP:
\[ \quad (x \circ y)^{-1} = x^{-1} \circ y^{-1} \mbox{ for all } x,y \in X. \]

This identity implies, that the inverse map $\iota: X \to X, x \mapsto x^{-1}$ is an automorphism of the loop. 
Bol loops with AIP are called {\em  Bruck loops} and generalise abelian groups. In the literature Bruck loops
occur also under other names, such as $K$-loops [Kreuz, Kiech] or gyrocommutative gyrogroups [Ung].

Glauberman showed in \cite{GG1} and \cite{GG2}, that Bruck loops of odd order behave very well: 
These loops are soluble and allow generalisations of many theorems of group theory.
His famous $Z^\ast$-theorem was originally a byproduct of this work.

Then, forty years later, Aschbacher, Kinyon and Phillips  showed that the following holds in finite
Bruck loops \cite{AKP}: 
\begin{itemize}
\item Elements of 2-power order and elements of odd order commute in a more general sense, see Theorem~2 of [AKP].
\item Bruck loops are a central product of a subloop of odd order and a subloop generated by elements of 2-power order. 
\item Simple Bruck loops are of 2-power exponent.
\item The structure of minimal simple Bruck loops is very restricted.
\end{itemize}

This leaves the Bruck loops of 2-power exponent to be studied.
Notice, that in the simplest case of Bruck loops of exponent 2, the automorphic inverse property is already a
consequence of the Bol identity and the exponent 2 assumption.

In \cite{A} Aschbacher gave powerful restrictions on the structure of minimal simple Bol loops of exponent 2. 
Using the restrictions given in Aschbacher's paper, Nagy and independently Baumei\-ster and Stein 
found a simple Bol loop of exponent 2 and size 96 in April 2007 [Nag1, BS].
Furthermore Nagy produced an infinite sequence of simple Bol loops of exponent 2 [Nag1].
A bit later a simple Bruck loop of exponent 4 and size 96 was found by Baumeister and Stein [BS].  

Thus the weakening of the associativity law produces lots of generalised elementary abelian groups, which are simple Bruck loops and live in non-soluble groups.
Recall, that Aschbachers paper \cite{A} and its generalisations in \cite{AKP} restrict only the structure of minimal simple loops.
Suddenly the question arose, whether the class of Bruck loops is maybe as wild as the class of general Bol loops. 

In this paper we determine the structure of  a finite Bruck loop showing that the structure of a
finite Bruck loop is not as wild as suspected. The definition of a loop envelope and a twisted 
subgroup is given in the next section. Recall that for $G$ a group $O(G)$ is the biggest normal subgroup of $G$
of odd order and that $O^{2^\prime}(G)$ is the smallest normal subgroup of $G$ such that $G/O^{2^\prime}(G)$ is 
a group of odd order.

\begin{thm}
\label{direct product}
Let $X$ be a finite Bruck loop. 
Then the following holds.
\begin{itemize}
\item[(a)] $X = Y \times Z$ where $Y$ is a subloop with  $|Y|$ odd and $Z$ a subloop of 2-power exponent.
\item[(b)] Let $G$ be the enveloping group of $X$. Then $G = O(G) \times O^{2^\prime}(G)$.
\item[(c)]
A loop envelope $(G,H,K)$ of $Z$ where $H$ acts on $K$ and $K$ is a twisted subgroup consisting of
$2$-power elements satisfies the following.
\begin{itemize}
\item[(1)] $\overline{G} =G/O_2(G) \cong D_1 \times D_2 \times \cdots \times D_e$ with $D_i \cong \PGL_2(q_i)$ for $q_i\ge 5$ a Fermat prime 
or $q_i=9$ and $e$ a non-negative integer,
\item[(2)]
$D_i \cap \overline{H}$ is a Borel subgroup in $D_i$,
\item[(3)] $F^\ast(G) = O_2(G)$,
\item[(4)] $ \overline{K}$ is the set of involutions in $\overline{G}\setminus{\overline{G}^\prime}$.
\end{itemize}
\end{itemize}
\end{thm}

\begin{rem}
\begin{itemize}
\item[(1)] Notice, that it may be $e = 0$ in which case $G$ is a $2$-group and $X$ a soluble loop!
\item[(2)] If $(G,H,K)$ is the Baer envelope of $Z$, then it satisfies the assumptions required in the
theorem.
\end{itemize}
\end{rem}

To prove Theorem~\ref{direct product} we also show that if 
$X$ is a finite Bruck loop with loop envelope $(G,H,K)$, then $O_2(G)$ is the group to a subloop of $X$.

A direct consequence of the theorem is the following.

\begin{coro}\label{characterisation soluble}
Let $X$ be a finite Bruck loop with enveloping group $G$. Then $X$ is soluble if and only if
$G= O(G) \times O_2(G)$.
\end{coro}

Now not only the class of Bruck loops is understood better,
but also the more general class of Bol-$A_r$-loops.

\begin{coro}\label{Ar-loops}
 Let $X$ be a Bol-$A_r$-loop. Then there is a normal subloop $Y$ of $X$ which is a group 
such that $X/Y$ is as described in Theorem~\ref{direct product}.
\end{coro}

As a group theoretic corollary we obtain:

\begin{coro}
Let $G$ be a finite group and $H \le G$, such that there is a transversal $K$ to $H$ in $G$ which
is the union of $1 \in G$ and $G$-conjugacy classes of involutions. If $G= \langle K \rangle$, then $(G,H,K)$
is a loop envelope to a Bruck loop of  exponent $2$ with $H$ acting on $K$. Therefore, Theorem~\ref{direct product}(c) describes $G$, $H$ and $K$.
\end{coro}

Moreover, we show Sylow's Theorem for the prime $2$.

\begin{thm}\label{Sylow2} [Sylow's Theorem]
Let $X$ be a finite Bruck loop. 
\begin{itemize}
 \item[(1)]  There is a subloop $P$ of $X$ such that $|P|=|X|_2$.
 \item[(2)] All subloops of $X$ of size $|X|_2$ are conjugate under $H$, the group of inner automorphisms of $X$.
 \item[(3)]  If $Y\le X$ with $|Y|$ a power of $2$, then there is an $h \in H$ such that $Y \le P^h$. 
\end{itemize}
\end{thm}

\begin{rem}
In fact, Theorems 12 and 14 of [Glaub2] as well as Theorem~\ref{direct product} yield if $p$ is an odd prime 
which divides the order of $X$ and which does not divide $q+1$ for any 
Fermat prime $q$ or $q = 9$, then there is a subloop $P$ of $X$ such that $|P|=|X|_p$.
\end{rem}

We also get Lagrange's Theorem for Bruck loops.

\begin{thm}\label{Lagrange} [Lagrange's Theorem]
Let $X$ be a finite Bruck loop and $Y \le X$ a subloop. 
Then $|Y|$ divides $|X|$. 
\end{thm}

Finally, the Theorem of Hall holds as well:

\begin{thm}\label{Hall} [Hall's Theorem]
Let $X$ be a finite Bruck loop and let $\Pi$ be the set of primes dividing the order of $X$.
 Then $X$ is soluble if and only if there is a Hall $\pi$-subloop in $X$ for every subset
$\pi$ of $\Pi$.
\end{thm}

There is even a stronger version of that theorem:

\begin{thm}\label{Hall2} 
 Let $X$ be a finite Bruck loop.
 Then $X$ is soluble if and only if there is a Sylow subloop in $X$ for every prime dividing
$|X|$.
\end{thm}

The organisation of the paper is as follows. In the next section we recall the relevant definitions
and the notation. Then in the third section we collect our previous results which will be needed in the proofs of the theorems. In Section~4 we prepare the proof of Sylow's Theorem for the special case of Bruck loops of $2$-power exponent by  calculating the number of elements in the intersection of $K$ with a Sylow $p$-subgroup. 
Moreover, we prove Theorems~\ref{Sylow2} and \ref{Hall} in the special case of Bruck loops of 
$2$-power exponent in that section. 
Finally 
Theorems~\ref{direct product}, ..., \ref{Hall2} will be shown in full generality in the last section.

\section{Definitions and Notation}

We follow the notation of Aschbacher \cite{A} and \cite{AKP}. Please keep in mind that in this paper
we study exclusively finite loops.

Baer observed that loops can be translated into the language of group theory \cite{B}.
This translation is as follows.
Let $X$ be a loop and let 
$$\rho: X \to Sym(X),~~~x \to \rho_x.$$
Define 
$$G := \RMult(X) := \langle \rho_x~|~ x \in X\rangle \le Sym(X)$$
$$H:=\Stab_G(1),~\mbox{where}~1 \in X,~\mbox{and}~K:=\{ \rho_x~|~ x \in X \}.$$ 

Then 

\begin{enumerate}
\item[(1)] $1 \in K$ and $K$ is a transversal to all conjugates of $H$ in $G$.
\item[(2)] $H$ is core free.
\item[(3)] $G = \langle K \rangle $.
\end{enumerate}

The group $G$ is called the {\em enveloping group of $X$} (or right multiplication group) and the triple $(G,H,K)$  the
 {\em Baer envelope of the loop}.
Baer also observed that whenever  $(G,H,K)$ is a triple with $G$ a group, $H \le G$ and $K \subseteq G$
satisfying condition (1),
then we get a loop on $K$ by setting $x \circ y = z$, $x,y \in K$ whenever $z$
is the element in $K$ such that $Hxy = Hz$.  This loop is called the {\em loop related to
$(G,H,K)$}.

The triple $(G,H,K)$ with $G$ a group, $H \le G$ and $K \subseteq G$ is called 
a {\em loop folder}, {\em faithful loop folder} or {\em loop envelope}
if (1), (1) and (2) or (1) and (3) hold, respectively.
In general there are  many different loop folders to a given loop.

If $X$ is a Bol loop and $(G,H,K)$ the Baer envelope of $X$, then $K$  is a {\em twisted subgroup},
that is $1 \in K$ and whenever $x,y \in K$, then $x^{-1}$ and $xyx$ is in $K$.
If, moreover, $X$ is a Bruck loop, then $H$ acts on $K$ by conjugation, [Asch, 4.1].
A {\em Bruck folder} is a loop folder $(G,H,K)$, if the following holds

\begin{itemize}
\item[(1)] $K$ is a twisted subgroup	
\item[(2)] $H$ acts on $K$ by conjugation.
\end{itemize}

We say that a Bruck folder is a {\em BX2P-folder}, if also

\begin{itemize}
\item[(3)] The elements in $K$ are of $2$-power order
\end{itemize}

If $(G,H,K)$ is a {\em BX2P-folder}, then the loop related to it is a Bruck loop of $2$-power exponent [BSS].
Moreover, notice that the subgroup $H$ of the Baer envelope induces automorphisms on $X$ in a Bruck loop. These are
called the {\em inner automorphisms of $X$}.

Subloops,  homomorphisms, normal subloops, factor loops and simple loops are defined as usual in universal algebra:
A {\em subloop} is a  nonempty subset which is closed under loop multiplication.

{\em Homomorphisms} are maps between loops which commute with loop multiplication. The map defines an equivalence relation
on the loop, such that the product of equivalence classes is again an equivalence class. 
{\em Normal subloops} are preimages of 1 under a homomorphism and therefore subloops.
A normal subloop defines a partition of the loop into blocks (cosets), such that the set of products of elements from two blocks is again a block.
Such a construction gives factor loops as homomorphic images with the block containing 1 as the kernel. 
{\em Simple loops} have only the full loop and the 1-loop as normal subloops.

For instance if $(G,H,K)$ is a loop folder defining a loop $X$ and $G_0$ a normal subgroup of $G$ which contains $H$,
then $(G_0, H, G_0 \cap K)$ is a loop folder to a normal subloop $X_0$ of $X$.

A {\em subfolder} $(U,V,W)$ is a loop folder with $U \le G$, $V \le U \cap H$ and $W \subseteq U \cap K$.
The loop to a subfolder of a loop folder is a subloop  of the loop to the loop folder.

Finally we recall the definition of a soluble loop given in [Asch]. A loop $X$ is {\em soluble} if there
exists a series $1 = X_0 \leq \cdots \leq X_n= X$ of subloops with $X_i$ normal in $X_{i+1}$ and
$X_{i+1}/X_i$ an abelian group.

Let $\pi$ be a set of primes. A natural number $n$ is a {\em $\pi$-number} if $n =1$ or $n$ is the product of
powers of primes in $\pi$. Assume that $X$ is a loop such that every element of the loop generates a group.
We say that $X$ is a $\pi$-loop, if the order of $X$ is a $\pi$-number. Notice that this definition is 
different from the one given in [Glaub1].  For loops of odd order these two concepts coincides (see
[Glaub1, p. 394, Corollary 2]), but not for loops of even order (see the Aschbacher loop in [BS]).

In order to distinguish the two concepts we propose to use the following notations:
A {\em local $\pi$-loop} is a loop such that the orders of the elements are all $\pi$-numbers
and a {\em global $\pi$-loop} is a loop such that the order of the loop is a $\pi$-number.

Then there are local $2$-loops which are not global $2$-loops, see [BS].

We say that a subloop $Y$ of $X$ is a {\em $\pi$-Hall subloop}, if $|Y|_\pi = |X|_\pi$.

\section{Previous results}

In the following let $(G,H,K)$ be a loop folder. We can see from the structure of a subgroup $U$ of
$G$, if it defines a subloop.

\begin{lemma}
\label{subfolders} [BSS, 2.1, 2.2]
A subgroup $U \le G$ gives rise to a subfolder $(U,V,W)$, if and only if $U = (U \cap H)(U \cap K)$.
Then $V = U \cap H$ and $W = U \cap K$. 
In particular, subgroups of $G$ which contain either  $H$ or $\langle K \rangle$ give rise to subfolders.
\end{lemma}

\begin{lemma}
\label{BruckFolder}[BSS, 2.15 (2) - (4)]
Let $(G,H,K)$ be a Bruck folder. Then the following holds.
\begin{enumerate}
\item[(1)]  There exists a unique $\tau \in \Aut(G)$ with $[H,\tau]=1$ and $k^\tau = k^{-1}$ for all $k \in K$. 
\item[(2)]  The set $\Lambda= \tau K \subseteq \Aut(G)$ is $G$-invariant.
\item[(3)]  Subfolders and homomorphic images are Bruck folder.
\end{enumerate}
\end{lemma}

Notice, that \ref{BruckFolder}(3) implies that subloops of Bruck loops are again Bruck loops.
In order to prove Sylow's Theorem we will work in a group slightly bigger than $G$.

\begin{defi}\label{Lambda}

Let $(G,H,K)$ be a Bruck folder and $\tau \in \Aut(G)$ the automorphism introduced in \ref{BruckFolder}(1). 
Then let \[{\bf  G^+} := G\langle \tau \rangle,\] the semidirect product of $G$ with $\tau$,
\[{\bf  H^+}:= H \langle \tau \rangle \le G^+~\mbox{and}~
{\bf  \Lambda }:= \tau K \subseteq G^+.\]  By \ref{BruckFolder}(1) and (2) $\Lambda$ is a $G^+$-invariant set of involutions.
\end{defi}

The following are powerful facts.

\begin{lemma}
\label{noHinvert}[BSS, 2.17, 3.3]
Let $(G,H,K)$ be a BX2P-folder and let $\overline{G}=G/O_2(G)$.
Then 
\begin{itemize}
\item[(1)] $k^2 \in O_2(G)$ for all $k$ in $K$. 
\item[(2)] $1 \in \overline{K}$ and $\overline{K}$ is a union of $\overline{G}$-conjugacy classes.
\item[(3)] Let $g \in G$ and $h \in H$.
If $(h^g)^k = (h^g)^{-1}$ for some $k$ in $K$, then $h^2 = 1$. 
\end{itemize}
\end{lemma}

We already have some information on soluble Bruck loops, see also [AKP, Corollary 4].

\begin{lemma}
\label{soluble_loops} [BSS, 3.8, 3.9, 3.10]
Let $(G,H,K)$ be a  BX2P-envelope to a Bruck loop $X$ of $2$-power exponent. Then the 
following holds
\begin{itemize}
 \item[(1)] $X$ is soluble if and only if $|X|$ is a power of $2$.
 \item[(2)] If $X$ is soluble, then $G$ is a $2$-group.
 \item[(3)] If $G = O_2(G)H$, then $X$ is soluble.
\end{itemize}
\end{lemma}

In [BSS] we introduced the concept of passive groups.
\medskip\\
\noindent
{\bf Definition}
 A finite nonabelian simple group $S$ is called {\em passive}, if whenever $(G,H,K)$ is a BX2P-folder with 
\[ F^\ast(G/O_2(G)) \cong S, \]
then $G = O_2(G) H$. 
\medskip\\

In that case consequently the loop to $(G,H,K)$ is of 2-power size and  soluble  by \ref{soluble_loops}(3).
The main theorems of [BSS], respectively [S] are the following.

\begin{thm}
\label{EnvelopeGroups} [BSS, Theorem 1]
Let $(G,H,K)$ be a loop envelope to a Bruck loop of exponent $2$ such that $K$ is a twisted subgroup, such that $H$ acts on $K$ and such that $K$ consists of elements of $2$-power order.
Moreover, assume that 
every non-abelian simple section of $G$ is either passive or isomorphic to $\PSL_2(q)$
 for $q=9$ or $q \ge 5$ is a Fermat prime.
Then the following holds
\begin{itemize}
\item[(1)] $\overline{G}:= G/O_2(G) \cong D_1 \times D_2 \times \cdots \times D_e$ with $D_i \cong \PGL_2(q_i)$ for $q_i\ge 5$ a Fermat prime 
or $q_i=9$ and $e$ a non-negative integer,
\item[(2)]
$D_i \cap \overline{H}$ is a Borel subgroup in $D_i$,
\item[(3)] $F^\ast(G) = O_2(G)$,
\end{itemize}
\end{thm}

\begin{thm}\label{MainClassificationTheorem} [S, Theorem 1]
Let $S$ be a non-abelian simple group, then $S$ is passive or isomorphic to $\PSL_2(q)$
 for $q=9$ or $q \ge 5$ a Fermat prime.
\end{thm}

It follows that the order of a Bruck loop of exponent $2$ is very restricted:

\begin{coro}\label{order}
Let $X$ be a Bruck loop of $2$-power exponent. Then \[|X| = 2^{a}\prod_{i=1}^e (q_i+1)\]
for some $e \in \NN\cup \{0\}$ and $q_i = 9$ or a Fermat prime. Moreover, $|X|_2 = 2^{a+e}$.
If $(G,H,K)$ is the Baer envelope of $X$, then $2^a = |O_2(G):O_2(G) \cap H|$. 
\end{coro}	

\section{Bruck loops of 2-power exponent}

In this section $(G,H,K)$ will always be a  BX2P-envelope to a non-soluble Bruck loop. 

By the definition of a BX2P-envelope it follows that $(G,H,K)$ is a loop envelope to a Bruck loop
of exponent $2$ such that 
$H$ acts on $K$, such that $K$ is a twisted subgroup and
such that the elements in $K$ are of $2$-power order. As by Theorem~\ref{MainClassificationTheorem}
every non-abelian simple non-passive group is isomorphic to $\PSL_2(q)$ for $q=9$ or $q \ge 5$ a Fermat prime,
Theorem~\ref{EnvelopeGroups} yields

\begin{lemma}
\begin{enumerate}
\item[(1)] $\overline{G}:= G/O_2(G) \cong D_1 \times D_2 \times ... \times D_e$ $~\mbox{with}~D_i \cong \PGL_2(q_i)$
and $q_i=9$ or a Fermat prime $q_i \ge 5$. 
\item[(2)] $D_i \cap \overline{H}=:B_i$ is a Borel subgroup of $D_i$.
\end{enumerate}
\end{lemma}

Let $\pi_i$ be the projection of $\overline{G}$ onto $D_i$. 
We first study $\overline{H}$ and $\overline{K}$ in more detail.

\begin{lemma}
\label{HisBorel}
\begin{enumerate}
\item[(1)] $\overline{H}= \prod\limits_{i=1}^e B_i$.
\item[(2)] If $\overline{k} \in \overline{K}$ and $1 \le i \le e$, then $\pi_i(\overline{k})$ is either $1$ or an involution in 
$D_i\setminus{D_i^\prime}$. 
\end{enumerate}
\end{lemma}

\begin{bew}
By Theorem \ref{EnvelopeGroups} $$B:=\prod\limits_{i=1}^e B_i \le \overline{H}.$$ 
The involutions in $D_i^\prime$ invert elements of odd order which are conjugate to elements in $B_i$. Therefore, 
\ref{noHinvert}(2) and (3) imply (2).

Next we aim to show $\overline{H}= B$.
As $B_i$ is a maximal subgroup of $D_i$, it follows $\pi_i(\overline{H}) = B_i$ or $D_i$, for $1 \le i \le e$.

Assume that $B < \overline{H}$. Then $\pi_j(\overline{H}) = D_j$ for some $1 \le j \le e$. This implies
 $\langle B_j^{\overline{H}} \rangle = D_j$
is a subgroup of $\overline{H}$.
Let $\overline{k}$ be an element in $\overline{K}$ and consider $\pi_j(\overline{k})$. If $\pi_j(\overline{k}) \ne 1$, then $\pi_j(\overline{k})$ inverts some element of odd prime order in $D_j$ by Baer-Suzuki. Then $\overline{k}$ inverts an element 
in $\overline{H}$ and therefore $k$ inverts an element in $H$ [Asch, (8.1)(1)].
This yields  a contradiction to \ref{noHinvert}. Therefore, $\pi_j(\overline{k})=1$ for all $\overline{k} \in \overline{K}$, 
which contradicts $\overline{G}= \langle \overline{K} \rangle$.

This shows $B = \overline{H}$.
\end{bew}

\subsection{Some subloops of Bruck loops of $2$-power exponent}

Now we can prove, that $O_2(G)$ is a group to a subloop:

\begin{lemma}
\label{O2subloop}
$O_2(G) H \cap K = O_2(G) \cap K$ and $O_2(G)= (O_2(G) \cap H)(O_2(G) \cap K)$.
\end{lemma}

\begin{bew}
By \ref{HisBorel}, $O_2(\overline{H})=1$. By \ref{subfolders} the subgroup $O_2(G)H$ defines a subloop, which is soluble by
 \ref{soluble_loops} (3).
Therefore $O_2(G)H = H(O_2(G)H \cap K)$ and 
 $\langle  O_2(G)H \cap K \rangle$ is a $2$-group by \ref{soluble_loops} (2), which yields 
$\langle  O_2(G)H  \cap K \rangle \le O_2(O_2(G)H) = O_2(G)$. Hence $O_2(G)H \cap K  = O_2(G)\cap K$
and the Dedekind identity implies the statement.
\end{bew}

Then application of Lemma~\ref{subfolders} shows that $O_2(G)$ is a group to a subloop.
There are lots of other subloops: Let $I:= \{1,2,...,e\}$ and  let $G_J$ be the full preimage of $\prod\limits_{j \in J} D_j$
for $J \subseteq I$.

\begin{lemma}
\label{manysubloops}
$G_J =(G_J \cap H)(G_J \cap K)$ for every $J \subseteq I$.
\end{lemma}

\begin{bew}
For $J=\emptyset$ this is \ref{O2subloop} and for $J=I$ this is the loop folder property. 

Let $x \in G_J$ and $x=h k$ with $h \in H$, $k \in K$. Let $l \in I-J$. As $\pi_l(x)=1$,
we cannot have $\pi_l(\overline{k}) \ne 1$: Else by \ref{HisBorel}(2), $\pi_l(\overline{k})$ is some involution of $\PGL_2(q_l)$ outside $\PSL_2(q_l)$.
But $\pi_l(\overline{H})=B_l$ and $B_l$ contains only involutions from $\PSL_2(q_l)$. 
So $\pi_l(\overline{k})=1$, thus $\pi_l(\overline{h})=1$ too. This implies the statement.
\end{bew}

\subsection{Preparations for Sylow's Theorem}

Our next goal is to produce subloops to certain Sylow-2-subgroups $P$ of $G$. Therefore we calculate $|P^+ \cap \Lambda|$. 

\begin{lemma}
For every $J \subseteq I$, $\overline{G}$ has a unique conjugacy class ${\cal C}_J$ of elements such that whenever
$t \in {\cal C}_J$, then  $\pi_i(t)=1$ for $i \not\in J$ and $\pi_i(t)$ is some involution in $D_i\setminus{D_i^\prime}$ 
for $i \in J$.
Moreover \[ |{\cal C}_J|= \prod_{j \in J} q_j \frac{q_j-1}{2}. \] 
\end{lemma}

\begin{bew}
This is immediate from the structure of $\overline{G}$. Recall, that for $q$ odd, the centraliser of an involution in $\PGL_2(q)$ is
the normaliser of a torus of size either $q-1$ or $q+1$. In our case $q-1$ is divisible by 4, so inner involutions of $\PSL_2(q)$
have a centraliser of size $2(q-1)$ while outer involutions have centraliser size $2(q+1)$.
\end{bew}

For $J \subseteq I$  let $t \in {\cal C}_J$. We denote by $O_2(G^+)t$ the full preimage of $t$ in $G^+$. The number
$n_J:= |O_2(G^+)t \cap \Lambda|$ is well defined and independent of the choice of $t\in {\cal C}_J$. 
Then $$n_\emptyset  = |O_2(G^+) \cap \Lambda| = |O_2(G) \cap K|=  |O_2(G):O_2(G) \cap H|$$ by \ref{manysubloops}.

\begin{lemma}
\[ n_J = \frac{n_\emptyset \cdot 2^{|J|}}{ \prod\limits_{j \in J} (q_j-1) } \]
\end{lemma}
\begin{bew}
By \ref{manysubloops} $G_J$ defines a subloop, so $|G_J: G_J \cap H|=|G_J \cap K| = |G_J^+ \cap \Lambda|$.
As $$|G_J:G_J \cap H| = |\overline{G_J}:\overline{G_J} \cap \overline{H}| |O_2(G) :O_2(G) \cap H|,$$ we have
\[ |G_J:G_J \cap H|= n_\emptyset \prod_{j \in J} (q_j+1). \]

On the other hand 
\[ |G_J^+ \cap \Lambda|= \sum_{L \subseteq J} n_L |{\cal C}_L|. \]
We therefore get a system of equations for the $n_J$. 

Now the statement can be shown by induction on $|J|$.
For example for $|J|=1$ we get the equation $n_\emptyset (q_j+1) = n_\emptyset + n_{\{j\}} \cdot q_j \frac{q_j-1}{2}$,
which gives $n_{\{j\}} = \frac{2 n_\emptyset} {q_j-1} $. 
In general we have:
\[ n_\emptyset \prod_{j \in J} (q_j+1) = \sum_{L \subseteq J} n_L \prod_{j \in L} q_j \frac{q_j-1}{2}. \]
For $L \subseteq J$, $L \ne J$ we have the formula for $n_L$ by induction. On the other hand for any
numbers $q_j, j \in J$ the equation 
\[ \prod_{j \in J} (q_j +1) = \sum_{L \subseteq J} \prod_{j \in L} q_j \]
holds. After some calculation this gives exactly the formula for $n_J$.
\end{bew}

\begin{lemma}
\label{SylowLoopExists}
Let $P \in \Syl_2(G)$. Then $$|P \cap K|=|P^+ \cap \Lambda|= 2^e n_\emptyset = |G:H|_2=|X|_2.$$
If $P \cap O_2(G)H \in \Syl_2(O_2(G)H)$, then $P= (P \cap H)(P \cap K)$.
\end{lemma}

\begin{bew}
Notice that $|P^+ \cap \Lambda|$
is independent of the choice of $P$, as $\Lambda$ is $G^+$-invariant. 

Furthermore, $|P \cap K| = |P^+ \cap \Lambda|$ as $\tau \in O_2(G^+) \le P^+$.

We choose $P \in \Syl_2(G)$ with $P \cap O_2(G)H \in \Syl_2(O_2(G)H)$. 
Then also $P^+ \cap O_2(G^+)H^+ \in \Syl_2(O_2(G^+)H^+)$.

Let $i \in I$ and consider $P_i=\pi_i(\overline{P}) \in \Syl_2(D_i)$. Then $P_i$ is a dihedral group,
$P_i \cap \overline{H}$ is a cyclic group of size $q_i-1$. The other coset of $P_i \cap \overline{H}$
in $P_i$ consists entirely of involutions, half of them involutions in $D_i^\prime$ and half of them in $D_i\setminus{D_i^\prime}$. 
As all involutions in $D_i\setminus{D_i^\prime}$ are conjugate in $\overline{G}$, it follows
\[ \pi_i(\overline{P^+}) \cap \overline{\Lambda} = 1 + \frac{q_i-1}{2}, \]
where $1$ is a summand as $1 \in \overline{\Lambda}$. 
This shows for $J \subseteq I$: \[ |\overline{P^+} \cap {\cal C}_J| = \prod_{j \in J} \frac{q_j-1}{2}. \]
As 
\[ |P^+ \cap \Lambda| = \sum_{J \subseteq I} n_J |\overline{P^+} \cap {\cal C}_J|, \]
it follows that
\[ |P^+ \cap \Lambda| = \sum_{J \subseteq I} \frac{n_\emptyset 2^{|J|}}{\prod\limits_{j \in J} (q_j-1)} \prod_{j \in J} \frac{q_j-1}{2} = 2^{|I|} n_\emptyset = 2^e n_\emptyset. \]
By the Dedekind identity we have $O_2(G^+) (P^+ \cap H^+)= P^+ \cap O_2(G^+)H^+$. This gives
\[ \frac{|O_2(G^+)||P^+\cap H^+|}{|O_2(G^+) \cap P^+ \cap H^+|}= |P^+\cap O_2(G^+)H^+|. \]
Theorem~\ref{EnvelopeGroups} yields
$$|P^+ \cap O_2(G^+)H^+| = \frac{|G^+|_2}{|G^+:O_2(G^+)H^+|_2} = \frac{|G^+|_2}{2^e}.$$ 
As $O_2(G) \leq P$ and 
$$|O_2(G)|/|O_2(G) \cap P \cap H|
= |O_2(G):O_2(G) \cap H| =n_\emptyset,$$
 it follows that
$$|P^+ \cap H^+|=|P \cap H|=\frac{|G|_2}{2^e n_\emptyset}$$
 and therefore,
$$|P:P \cap H|= 2^e n_\emptyset.$$
Hence $|P:P \cap H| =  |P \cap K|$, which yields $P= (P \cap H)(P \cap K)$.

Moreover, $$|X|_2=|K|_2 = |G:H|_2 = |G:O_2(G)H|_2 |O_2(G)H:H|_2$$
$$ = 2^e |O_2(G):O_2(G) \cap H| = 2^e n_\emptyset.$$
\end{bew}

\subsection{Sylow's theorem for Bruck loops of $2$-power exponent}

Next we show, that the subloops of size $|X|_2$ have some nice properties. Therefore, we need to recall
some facts about $\PGL_2(q)$.
 
\begin{lemma}
\label{PGLlemma}
Let $Z \cong \PGL_2(q)$ with $q=9$ or $q\ge 5$ a Fermat prime. Let $B$ be a Borel subgroup of $G$ and ${\cal C}$
the class of involutions in $Z\setminus{Z^\prime}$. 
\begin{enumerate}
\item[(1)] $B$ has two orbits on $\Syl_2(Z)$: one orbit of size $q$ and one of size $\frac{|B|}{2}$. 
\item[(2)] If $P \in \Syl_2(Z)$, then either $P \cap B \in \Syl_2(B)$ or $|P \cap B|=2$.
\item[(3)] Let $A \subseteq \{1 \} \cup {\cal C}$ and suppose that 
$D= \langle A \rangle$ is a $2$-group such that $D=(D \cap B)A$.
Then there is a $Q \in \Syl_2(Z)$ such that $D \le Q$ and $Q \cap B \in \Syl_2(B)$. 
\end{enumerate}
\end{lemma}

\begin{bew}
Let $\Omega$ be the set of points of the projective line related to $Z$. Then $|\Omega| = q+1$ and 
$Z$ acts triply transitive on $\Omega$.
Moreover, $B$ is the stabiliser in $Z$ of a point $a$
of $\Omega$ and every Sylow $2$-subgroup of $G$ is the setwise stabiliser of two points of $\Omega$.

It follows that $B$ has two orbits on the set of pairs: one consisting of the pairs containing $a$ and the other one
consisting of those not containing $a$. Their length are $q$ and $q(q-1)/2$, respectively. This shows (1).

If $P \in \Syl_2(Z)$ fixes a pair in the first orbit, then $P \cap B \in \Syl_2(B)$.  If $P$ fixes a pair in the second
orbit, then $P \cap B$ fixes a point and a pair of points setwise and is therefore just an involution, which is (2).

As $D=(D \cap B)A$, it follows that $|D:D\cap B| = 2$. Thus $D\cap B$ fixes the point $a$ and $a^D$ is of length $2$.
Let $Q$ be the stabiliser of $a^D$ in $G$. Then $Q \in \Syl_2(Z)$, $D \le Q$ and $Q \cap B \in \Syl_2(B)$, which is (3).
\end{bew}

The following is fundamental for the proof of the $2$-Sylow Theorem.

\begin{lemma}
\label{goodSylow}
Let  $(G,H,K)$ be a faithful BX2P-envelope and $U$ a $2$-subgroup of $G$ such that 
\begin{itemize}
\item $U = \langle U \cap K \rangle$  
\item $U=(U \cap H)(U \cap K)$.
\end{itemize}
Then there is a Sylow-2-subgroup $Q$ of $G$ such that $U \le Q$ and $Q \cap O_2(G) H \in \Syl_2(O_2(G)H)$. 
\end{lemma}

\begin{bew}
For fixed $1 \leq i \leq e$ let $$Z:= \pi_i(\overline{G}), D:=\pi_i(\overline{U}), B:=\pi_i(\overline{H})~\mbox{and}~A:=\pi_i(\overline{U} \cap \overline{K}).$$
 By \ref{HisBorel} $${\cal C}:=\pi_i(\overline{K}) -\{1 \}$$ is the class of involutions in 
$Z\setminus{Z^\prime}$. Moreover by the homomorphism property of $\pi_i$ it follows that 
$$\pi_i(\overline{U}) = \pi_i(\overline{U} \cap \overline{H}) \pi_i(\overline{U} \cap \overline{K}).$$ 
This yields, as $$\pi_i(\overline{U} \cap \overline{H}) \le \pi_i(\overline{U}) \cap \pi_i(\overline{H})= D \cap B,$$
that $D = (B\cap B)A$.
Hence,  $Z,B,A,D$ satisfy the assumptions of \ref{PGLlemma}(3).
Therefore, \ref{PGLlemma}(3) implies that $\pi_i(U)$ is contained in a Sylow $2$-subgroup $Q_i$ of 
$\pi_i(\overline{G})$ and that
$Q_i \cap \overline{H}$ is a Sylow $2$-subgroup of  $\pi_i(\overline{H})$. 

Let $Q$ be the preimage of $\prod\limits_{i \in I}Q_i$. Then $U \leq Q$ and $$Q \cap O_2(G) H \in \Syl_2( O_2(G) H)$$
as asserted.
\end{bew}

\begin{coro}
\label{SylowLoops}
Let $X$ be a finite Bruck loop of $2$-power exponent and 
$Y$ a soluble subloop of $X$. 
\begin{enumerate}
\item[(1)] Then there is a subloop $Z$ of $X$ such that $Y \le Z$ and $|Z|=|X|_2$.
\item[(2)] All subloops of $X$ of size $|X|_2$ are conjugate under $H$.  
\end{enumerate}
\end{coro}
\begin{bew}
 If $Y$ is soluble, then $Y$ is a $2$-loop by \ref{soluble_loops}(1). 
Let $(G,H,K)$ be a faithful BX2P-envelope to $X$. Then there is a subgroup $U$ of $G$ such 
that 
\begin{itemize}
\item $U = \langle U \cap K \rangle$
\item $U = (U \cap H)(U \cap K)$ by \ref{subfolders}
\item $U$ is a $2$-group by \ref{soluble_loops}(2)
\end{itemize}
Hence by \ref{goodSylow} there is a Sylow $2$-subgroup $Q$ of $G$ such that $U$ is  contained in $Q$
and such that $Q \cap O_2(G)H \in Syl_2(O_2(G)H)$. Now \ref{SylowLoopExists} and \ref{subfolders} imply that 
$(Q,Q \cap H, Q \cap K)$ is a subfolder of our chosen folder. Let $Z$ the subloop of $X$ related
to that subfolder. As $Q \cap O_2(G)H$ is a Sylow $2$-subgroup of $O_2(G)H$, the intersection $Q \cap H$ is
a Sylow $2$-subgroup of $H$. As $Q = (Q \cap H)(Q \cap K)$ it follows that $|Q \cap K| = |X|_2$, which
proves (1).

Let $Y_2$ be a subloop of $X$ of size $|X|_2$. Then $Y_2$ is soluble and therefore by (1) there is Sylow $2$-subgroup $P$ of $G$
such that $(P, P \cap H, P\cap K)$ is a subfolder to a subloop $Z_2$ which contains $Y_2$ and which is of order $|X|_2$.
This shows that $Y_2 = Z_2$. Recall that $G/O_2(G) = D_1 \times \cdots \times D_e$ with $D_i \cong \PGL_2(q_i)$ and
$e \geq 0$. Then, as $P \cap H$ is a Sylow $2$-subgroup of $H$ (see also \ref{PGLlemma}(2)),  according to \ref{PGLlemma}(1) 
there is an element $h$ in $H$ which maps $Q$ onto $P$.
\end{bew}

Now we can easily prove Theorem~\ref{Sylow2}.
\bigskip
\\
\noindent
{\bf Proof of Theorem~\ref{Sylow2} for Bruck loops of $2$-power exponent.}
As the set consisting of the $1$-element of $X$ is a soluble subloop of $X$, Corollary~\ref{SylowLoops} yields (1).
(2) is the second statement of the corollary.

If $Y$ is a subloop of $X$ of $2$-power order, then $Y$ is soluble by \ref{soluble_loops}(1). Therefore (3) follows
from Corollary~\ref{SylowLoops} as well.
\hfill\qed

\subsection{Lagrange's theorem}

Now we can prove Lagrange's theorem for Bruck loops of $2$-power exponent.
\bigskip
\\
\noindent
{\bf Proof of Theorem~\ref{Lagrange} for $X$ a Bruck loop of $2$-power exponent.}
First of all notice, that by \ref{BruckFolder}(3) $Y$ is a Bruck loop of $2$-power exponent.
By \ref{SylowLoopExists}, we have $|Y|_2 \le |X|_2$: There is  a subloop of $Y$ of size $|Y|_2$ by 
Theorem~\ref{Sylow2} in case of Bruck loops of $2$-power exponent,
which is soluble by \ref{soluble_loops}.
Let  $U\le G$ be the $2$-group related to this subloop; so $|U \cap K| = |Y|_2$. As
by \ref{SylowLoopExists} $|P \cap K| = |X|_2$ for any Sylow-2-subgroup of $G$,
$|Y|_2$ is a divisor of $|X|_2$. 

Suppose $Y$ is nonsoluble. Then $|Y|_{2'} \ne 1$. There is a subgroup
$U\le G$ such that $U =(U \cap H)(U \cap K)$, $U = \langle U \cap K \rangle$  and $|Y|=|U:U \cap H| = |U \cap K|$. 
By \ref{soluble_loops} $U$ is not a soluble group. 
In the following we apply Theorem~\ref{EnvelopeGroups} on $U$ and on $G$ and use the notation introduced there. The map $\theta: U \to G:  u  \mapsto O_2(G) u$ gives a homomorphism from $U$ into $\overline{G}$ and
an injection from $U/(O_2(U)\cap O_2(G))$ into $\overline{G}$. 

Assume there is $D_i \leq \overline{G}$ such that $\pi_i(\overline{U})$ is non-soluble and properly
contained in $D_i$. Then by Theorem~\ref{EnvelopeGroups} $\pi_i(\overline{U}) \cong \PGL_2(5)$ and $D_i \cong \PGL_2(9)$. Elements of odd order from $U \cap H$ map to elements of odd order in $\overline{H}$, which yields a contradiction
as $\pi_i(\overline{U\cap H}) \cong 5:4$ and $\overline{H}$ is a $\{2,3\}$-group.
Hence components of $U/O_2(U)$ project surjectively onto components of $G/O_2(G)$.
This implies together with \ref{order} that \[|Y| = 2^a\prod_{j\in J}(q_j+1)~\mbox{where}~J~\mbox{ is a subset of }~\{1, \ldots , e\}.\]
Thus the odd  part of $|Y|$ divides $|X|$. 

As we already saw that  $|Y|_2 \le |X|_2$, it follows that the order of $Y$ divides the order of $X$.
\hfill\qed

\section{The finite Bruck loops}

In this section we prove the main theorems.
\bigskip
\\
{\bf Proof of Theorem~\ref{direct product}.}
Let $X$ be a finite Bruck loop. Then according to [AKP, Theorem 1]
 \[X = O^{2^\prime}(X) * O(X),~\mbox{where}~Z:= O^{2^\prime}(X)\] 
is the subloop generated by all $2$-elements of $X$ and $$Y:= O(X)$$ the largest normal 
subloop of $X$ of odd order. Notice, that the definition of the subloop $Z$ is different from
that one in [AKP].
Let $(G,H,K)$ be the Baer-envelope of $X$ and set $$G_2:= O^{2^\prime}(G).$$
 Then $G_2 = \la \rho_x ~|~ x \in Z \ra \leq \Sym(X)$, see [AKP, 6.1].
Moreover, $G_1:= O^{2}(G)$ is the enveloping group of $O(X)$ and $G = G_1 *G_2$, see [AKP, Proof of 6.1].
Set $U:= G_1 \cap G_2$ and $T = Z \cap Y$. Then $U$ is a subgroup of $Z(G)$ and therefore, it acts
semiregularly on $Z$. 

Moreover, $Z/T$ is a Bruck loop of $2$-power exponent with enveloping group
$G_2/V$ where $V \leq U$ is the enveloping group of $T$, see [AKP] Theorem~1 and [Asch] 2.6. 
In particular, $V$ is a subgroup of $Z(G_2)$ of odd order.

Set $\tilde{G}_2=G_2/V$. Then
$\tilde{G}_2/O_2(\tilde{G}_2)  \cong D_1 \times \cdots \cong D_e$ with $D_i \cong \PGL_2(q_i)$, where
$q_i$ is a Fermat prime or $9$ by Theorems~\ref{EnvelopeGroups} and \ref{MainClassificationTheorem}

We claim that $G_2$ splits over $V$. Clearly, $O_2(G_2)V$ splits over $V$. Therefore, $G_2$ splits over $V$
if and only if $\overline{G}_2:=G_2/O_2(G_2)$ splits over $\overline{V}$. 

Moreover, $\overline{G}_2$ splits over $\overline{V}$ if and only if the full preimage $L_i$ of $D_i$ in $\overline{G}_2$ splits over
$\overline{V}$, for $1 \leq i \leq e$. As $L_i/\overline{V} \cong \PGL_2(q_i)$ and as $\overline{V}$ is of odd order it follows that the extension splits
or that $q_i = 9$ and $|L_i^\infty \cap \overline{V}| = 3$, see [Atlas, p. XVI, Table 5]. Assume the latter. Then every involution in 
$L_i\setminus{L_i^\prime}$ inverts $L_i^\infty \cap \overline{V}$, see the action of the automorphism $2_3$ of $\PSL_3(4)$
on the Schur-multiplier of order $3$ of $\PSL_3(4)$ in [Atlas, p. 23 ], which contradicts $V \leq Z(G_2)$. This shows that
$G_2$ splits over $V$.

Hence, as $G_2 = O^{2^\prime}(G)$, we get $V = 1$ as well as $T = 1$. Thus (a) and (b) follow.
The fact that $T =1$, also implies that $Z$ is a Bruck loop of $2$-power exponent. Part (c) is a direct consequence of
Theorem~\ref{EnvelopeGroups}.
\hfill\qed
\bigskip
\\
{\bf Proof of Theorem~\ref{Sylow2}.}
The assertion follows from Theorem~\ref{direct product} and the proof of Theorem~\ref{Sylow2} in the case 
of Bruck loops of $2$-power exponent.
\hfill \qed
\bigskip
\\
{\bf Proof of Theorem~\ref{Lagrange}.}
Following Bruck, it is enough to show, that this condition holds already in simple Bruck loops, see [Bruck, Chapter
V, p. 93, Lemma 2.1]. As simple Bruck loops
are either of prime order or of $2$-power exponent, we get the result from the proof of Theorem~\ref{Lagrange}
in the case of Bruck loops of $2$-power exponent.

A different proof of Theorem~\ref{Lagrange} without quoting Bruck is to apply  Theorem~\ref{direct product},
Corollary~4, p. 395, in [Glaub1]
 and 
 the proof of Theorem~\ref{Lagrange} in the case of Bruck loops of $2$-power exponent.
\hfill\qed
\bigskip
\\
{\bf Proof of Theorem~\ref{Hall}.}
Every finite Bruck loop of odd order is soluble, see Theorem 14 of [Glaub2], and contains therefore Hall $\pi$-subloops,
Theorem 12 [Glaub2]. If $X$ is a soluble Bruck loop, then by Theorem~\ref{direct product} $X$ is the direct product
of a Bruck loop of odd order and a Bruck loop of $2$-power order. Hence there is a $\pi$-subloop for every subset $\pi$ of
$\Pi$.

Now assume that $X$ is a Bruck loop such that there is a $\pi$-subloop for every subset $\pi$ of
$\Pi$.  If $X$ is of odd order, then it is soluble [Glaub2]. Assume now that $O^{2^\prime}(X) \neq 1$. If 
$O^{2^\prime}(X)$ is not of $2$-power order, then there is $q_i$, $q_i = 9$ or $q_i \geq 5$ a Fermat prime
such that a prime divisor $r\neq 2$ of $q_i+1$ divides $|O^{2^\prime}(X)|$, see Corollary~\ref{order}.
As there is no Bruck loop of $r$-power order, see Theorem~\ref{direct product}, it follows that
$O^{2^\prime}(X)$ is of $2$-power order and therefore soluble by \ref{soluble_loops}.
\hfill\qed
\bigskip
\\
{\bf Proof of Theorem~\ref{Hall2}.}
The previous proof also shows Theorem~\ref{Hall2}.
\hfill\qed

\subsection{Open questions}
There are still some open questions on Bruck loops:\\

\begin{itemize}
\item For which $q$ exist $M$-loops and/or $N$-loops as defined in \cite{AKP} and \cite{A}? 
There are only known examples for $q=5$.
\item Are there infinitely many Fermat primes ?
This is number theory...
\item What is the structure of simple Bruck loops in detail?
There are known examples with two composition factors of type $\Alt_5$. Are the nonabelian composition factors of $G$  pair wise isomorphic in a simple Bruck loop?
\item Is there a way to get the structure of $O_2(G)$ under control in $M$-loops, $N$-loops and/or Bruck loops of 2-power exponent
(For the definition of an $N$ and an $M$-loop see [BSS])?
\end{itemize}


\begin{thebibliography}{99}

\bibitem[Asch]{A} M.Aschbacher, On Bol loops of exponent 2, {\em J. Algebra} {\bf 288} (2005), 99-136

\bibitem[AKP]{AKP} M.Aschbacher, M.Kinyon, J.D.Phillips, Finite Bruck loops, {\em Transactions of the  AMS} {\bf 358}, No.7 (2006), 3061-3075 
%
\bibitem[ATLAS]{ATLAS} J.H.Conway, R.T.Curtis, S.P.Norton, R.A.Parker, R.A.Wilson {\em An ATLAS of finite groups} Oxford University Press, 1985
%
\bibitem[Baer]{B} R.Baer,  Nets and groups, {\em Trans. Amer. Math. Soc} {\bf 47} (1939), 110-141

\bibitem[BS]{BS} B.Baumeister, A.Stein, Self-invariant $1$-Factorizations of Complete Graphs and Finite Bol Loops of Exponent $2$, to appear in  {\em Beitr\"age zur Algebra und Geometrie}

\bibitem[BSS]{BSS} B.Baumeister, A.Stein, G.Stroth, On Bruck Loops of $2$-power Exponent, submitted 

\bibitem[Bol]{Bol} G.Bol {\em Gewebe und Gruppen} Math. Ann.  {\bf 114} (1937), 414-431
%
\bibitem[Bruck]{Bruck} R.H.Bruck, {\em A survey of Binary Systems}, Ergebnisse der Mathematik und ihrer Grenzgebiete.
Neue Folge, Heft 20. Reihe: Gruppentheorie, Springer Verlag, Berlin etc 1958
%

%
%
\bibitem[Glaub1]{GG1} G.Glauberman, On Loops of odd order, {\em J. Algebra} {\bf 1} (1964), 374-396 
%
\bibitem[Glaub2]{GG2} G.Glauberman, On Loops of odd order II, {\em J. Algebra} {\bf 8} (1968), 393-414
%
%
\bibitem[Lieb]{Liebeck} M.W.Liebeck, The classification of finite Moufang loops, {\em Math. Proc. Cambridge Philos. Soc} {\bf 102} (1987),33-47
%
\bibitem[Kiech]{Kiech} H. Kiechle, {\em The Theory of $K$-Loops},  Lecture Notes in Mathematics 1778, Springer, 
Berlin, Heidelberg, New-York, 2002.
%
\bibitem[Kreuz]{Kreuzer} A. Kreuzer, Inner mappings of Bruck loops, {\em Math. Proc. Cambridge Philos. Soc.}
{\bf 123} (1998), 53-57
%
\bibitem[Nag1]{nagy} G.Nagy, {\em A class of simple proper Bol loops}, {\em Preprint}
%
\bibitem[Nag2]{nagy-homepage} G.Nagy {\em Finite simple left Bol loops} \begin{verbatim}http://www.math.u-szeged.hu/~nagyg/pub/simple_bol_loops.html\end{verbatim}  
%
\bibitem[S]{Stein} A.Stein, On Bruck Loops of 2-power Exponent, II, Preprint 2009
%
\bibitem[Ung]{Ungar} A.A.Ungar, Beyond the Einstein Addition Law and its Gyroscopic Thomas Precession: The Theory of Gyrogroups and Gyrovectors Spaces Kluwer Academic Publishers, Doldrechts-Boston-London 2001  
 
\end{thebibliography}
\end{document}